\documentclass[11pt, reqno]{amsart}

 \usepackage{amsmath}
 \usepackage{amssymb}
\usepackage{bm}

\DeclareMathAccent{\mathring}{\mathalpha}{operators}{"17}

 \usepackage{color}

\newcommand{\mysection}[1]{\section{#1}
      \setcounter{equation}{0}}

\newtheorem{theorem}{Theorem}[section]
\newtheorem{lemma}[theorem]{Lemma}

\theoremstyle{definition}
\newtheorem{assumption}{Assumption}[section]
\newtheorem{definition}{Definition}[section]
\theoremstyle{remark}
\newtheorem{remark}{Remark}[section]

\newtheorem{example}{Example}[section]

\newcommand{\loc}{\text{\rm loc}}
\newcommand{\diam}{{\rm diam}}
\newcommand{\sign}{\text{\rm\,sign}\,}
\newcommand\inside{\rm int}

\newcommand\sfu{{\sf u}} 
\newcommand\sfv{{\sf v}}   

 \makeatletter 
 \def\dashint{%
 \operatorname%
 {\,\,\text{\bf--}\kern-.98em\DOTSI\intop\ilimits@\!\!}}
\def\ninf{\qopname\relax\@empty{inf\phantom{p}\!\!\!}}
 \makeatother

\newcommand\gb{\mathfrak{b}}

\newcommand\bbeta{\text{\raise-.2ex\hbox{$\bm{\beta}$}}}

\newcommand\bR{\mathbb{R}}

\def\bM{\mathbb{M}}

\newcommand\bS{\mathbb{S}}

\newcommand\cL{\mathcal{L}}

\newcommand\dist{{\rm dist}\,}

\newcommand\infsup{\operatornamewithlimits{inf\,\,\,sup}}

\begin{document}

\title[Elliptic equations with $L_{d}$-drift]
{Linear and fully nonlinear elliptic equations with $L_{d}$-drift}

\author{N.V. Krylov}
 
\email{nkrylov@umn.edu}
\address{127 Vincent Hall, University of Minnesota,
 Minneapolis, MN, 55455}

\keywords{Fully nonlinear equations, interior estimates,
solvability, unbounded coefficients}

\subjclass[2010]{35J60, 35J15}

\begin{abstract}
In subdomains of $\bR^{d}$ we consider uniformly elliptic equations
$H\big(v( x),D v( x),D^{2}v( x), x\big)=0$ with
the growth of $H$ with respect to $|Dv|$
controlled by the product of a function from $L_{d}$ times  $|Dv|$.
The dependence of $H$ on $x$ is assumed to be of BMO type.
Among other things we prove that there exists
$d_{0}\in(d/2,d)$ such that for any $p\in(d_{0},d)$
the equation with prescribed continuous boundary data
has a solution in class $W^{2}_{p,\loc}$.
Our results are new
even if $H$ is linear.
  
\end{abstract}

\maketitle

\mysection{Introduction and main results}

 In this article we consider elliptic  equations    
\begin{equation}
                                                \label{7.29.1}
H[v](x):= H\big(v( x),D v( x),D^{2}v( x), x\big)=0
\end{equation}
in subdomains $\Omega$ of $\bR^{d } $,   where $H(\sfu,x)$
is a function given for $x\in\bR^{d}$ and $\sfu=(\sfu',\sfu'')$,
$$
\sfu'=(\sfu'_{0},\sfu'_{1},...,\sfu'_{d}\big)\in\bR^{d+1},
\quad \sfu''\in\bS,
$$  
where $\bS$ is the set of symmetric $d\times d$-matrices.
 The ``coefficients''
of the first order derivatives of $v$  in \eqref{7.29.1}
are assumed to be in $L_{d}(\Omega)$ and  we take
$p\in(d_{0},d)$ for certain $d_{0}<d$.
We present some results 
about a priori estimates and the solvability in $W^{2}_{p,\loc}
(\Omega)$ of \eqref{7.29.1}. These results are new 
 even for linear equations (see Section \ref{section 10.10.2}
and Example \ref{example 1.21.1})
although in the linear case
results somewhat 
close to ours can be found in \cite{KK_19} 
under some additional regularity assumptions
on the matrix of second order coefficients
allowing one to rewrite the equation in divergence form.
Also see the references in \cite{KK_19}.
 Most likely our results are false  if $p=d$
even if the equation is linear.
 
In the literature, the   $W^2_{p,\loc}$, $p>d$, estimates 
like \eqref{1,25,3} with $\tau_{0}=0$
and $\Omega=B_{1}$
for viscosity solutions of a
class of fully nonlinear uniformly elliptic
equations of the form     
$$
H(D^2u,x)=f(x)
 $$
were first obtained by Caffarelli in \cite{Caf89} 
(see also \cite{CC_95}).  
His proof is
  based on an ingenious application of
 the Aleksandrov--Bakel'man--Pucci a
priori estimate, the Krylov--Safonov Harnack inequality, and
 a covering result which can be also found
in \cite{KS_80} and \cite{Sa_80}.
Our results are based on  
 ideas and results  
from \cite{Kr_13.1}, which uses the Evans-Krylov,
Fang-Hua Lin, and Fefferman--Stein theorems
as presented in \cite{Kr_18},
and results from recent papers   \cite{Kr_19_1},
and \cite{Kr_20}.
 By
exploiting a weak reverse H\"older's inequality, the result of
\cite{Caf89} was sharpened by Escauriaza in \cite{Es93}, 
who obtained the
interior $ W^2_p$-estimate for the same equations allowing
$ p>d-\varepsilon$, with a small constant
$\varepsilon>0$ depending only on the ellipticity constant and $d$. 
No terms with $Du$, however, were involved. In the present
article we use $p$ which is less than $d$
unlike  \cite{Kr_18}, where $p>d$ and the drift terms
are bounded.

The above cited works  \cite{Caf89} and
 \cite{CC_95}  are quite remarkable 
in one   respect--they do not suppose that $H$ is convex
 or concave in $ D^{2}u$  and relate
to any viscosity solution. 
The assumptions in \cite{Caf89} and \cite{CC_95}
are quite different from ours.
One of these assumptions is that the equations
$H(D^{2}u,x_{0})=0$ admit 
$C^{2}_{\loc}\big(B_{r}(x_{0})\big)$-solutions for any
$B_{r}(x_{0})\subset B_{1}$ and any continuous boundary data.   
Until now we only know that,  
generally, this assumption is satisfied 
if $H$ is convex or concave with respect to $\sfu''$.
 Paper \cite{SW_2016}
and the references there present a few exceptions. 

A number of existence results of $W^{2}_{p,\loc}$-solutions
and a priory estimates in $W^{2}_{p,\loc}$
obtained by means of the theory of viscosity solutions
can be found in \cite{Caf90}, \cite{CC_95}, and \cite{CCKS_96}.
In all of them $H$ is supposed to be Lipschitz continuous
in $\sfu$ uniformly with respect to $x$ and for any
$K>0$ and $|\sfu|\leq K$ to be
sufficiently
uniformly close to functions continuous 
 with respect to $x$. Note that these
assumptions exclude, for instance, Example \ref{example 1,13,1}
below,
and for, that matter, exclude {\em linear\/}
equations even with bounded coefficients and VMO-coefficients
in the main part.
On the other hand, our results do not cover those
from \cite{Caf90}, \cite{CC_95}, and \cite{CCKS_96}
either, in particular,
just because we are not dealing with viscosity solutions.
 
Note that in Theorem 4.2 of \cite{Wi09} one more
  interior estimate  of type \eqref{1,25,3} is obtained
under the assumptions that $H$ is convex in $\sfu''$, 
Lipschitz continuous in $\sfu$ and satisfies a 
continuity condition in $x$
similar to the one mentioned above. Again
some values of $p<d$ are allowed. Finally,
in \cite{CCKS_96} and \cite{Wi09} the function $H$
is assumed to be nonincreasing with respect
to $\sfu'_{0}$ unlike $H$ in our Theorem \ref{theorem 10.5.1}.
 
The article was motivated by Safonov's results
in \cite{Sa_10} where he, in particular, proved
the Harnack inequality and established
the H\"older continuity for harmonic functions
associated with linear elliptic equations with
measurable coefficients and drift in $L_{d}$.

To start the exposition of our results recall that
 $\bS$ is the set of symmetric $d\times d$ 
matrices and, for a fixed
$\delta\in(0,1]$,  let
$$
\bS_{\delta}=\big\{a\in\bS:\delta|\xi|^{2}\leq a^{ij}\xi^{i}\xi^{j}
\leq\delta^{-1}|\xi|^{2},\quad \forall \xi\in\bR^{d}\big\}.
$$

We fix a number $\|b\|<\infty$
and fix a nonnegative function $\gb\in L_{d}(\bR^{d})$
such that
$$
\|\gb\|_{L_{d}(\bR^{d})}\leq \|b\|.
$$

  Also fix some constants $K_{0},K_{F}\in[0,\infty)$
and fix a nonnegative $\bar G$ given on $\bR^{d}$.

Let $\Omega$ be an open bounded subset of $\bR^{d}$
satisfying the exterior ball condition.
Quite often we deal with
$$
\Omega^{\rho}=\big\{x\in\Omega:\dist(x,\partial\Omega)>\rho\big\},
$$
where $\rho>0$ is a given number.
For measurable $\Gamma\subset \bR^{d}$ we denote by $|\Gamma|$
the volume of $\Gamma$ 
 and if $f$ is a real-valued function on $\Gamma$
with finite integral,
then we set
$$
\dashint_{\Gamma}f\,dx=\frac{1}{|\Gamma|}\int_{\Gamma}f(x)\,dx.
$$

The following assumptions contain   parameters $\hat\theta,\theta\in(0,1]$
which are specified later in our results.
\begin{assumption}
                                    \label{assumption 10.5.1}
There are  Borel  
functions $F(\sfu,x)=F\big(\sfu'_{0},\sfu'',x\big)$  and $G(\sfu,x)$ 
   such that
$$
H =F +G .
$$
 Furthermore, for all $\sfu''\in\bS,\sfu'\in\bR^{d+1}$, and 
$x\in\bR^{d}$, we have 
\begin{equation}
                                                     \label{1,14,2}
\big|G(\sfu,x)\big|\leq \hat\theta|\sfu''|+ K_{0}|\sfu'_{0}|
+\gb(x)|[\sfu']|+\bar{G}(x),\quad 
F(0,x)\equiv0,
 \end{equation}
where
$$
[\sfu']:=(\sfu'_{1},...,\sfu'_{d}).
$$
 
\end{assumption}
 
Introduce
$$
B_{r}(x)=\big\{y\in\bR^{d}:|x-y|<r\big\},\quad B_{r}=B_{r}(0).
$$

Recall that   Lipschitz
continuous functions are almost everywhere differentiable,
thanks to the Rademacher theorem.

\begin{assumption}
                                \label{assump1} 
  (i) The function $F$ is Lipschitz continuous with respect to $\sfu''$
with Lipschitz constant $K_{F}$.

Moreover, there exist  $R_0\in(0,1]$ and $\tau_{0}\in[0,\infty)$
such that, if   
$r\in (0, R_0]$, $z\in\Omega$,  $B_{r}(z)\subset\Omega$, 
and $\sfu'_{0}\in\bR$, then
one can find a {\em convex\/} function $\bar{F} (\sfu'' )=
\bar{F}_{z,r,\sfu'_{0}} (\sfu'' )$ (independent
of $x $)  for which

(ii) We have $\bar{F}(0)=0$ and $ D_{\sfu'' }\bar{F} \in\bS_{\delta}$
at all points of differentiability
of $\bar{F}$;
 
(iii) 
For any $\sfu''\in\bS$ with $|\sfu''|=1$,  
 we have    
\begin{equation}
                                                \label{7.30.2}
\dashint_{ B_{r}(z)}\sup_{\tau>\tau_{0}}\tau^{-1}
\big|F\big(\sfu'_{0},\tau \sfu'' ,x\big)-\bar{F}(\tau \sfu'')\big| \,dx\leq \theta;
\end{equation}

(iv) There exists a continuous increasing 
function $\omega_{F }(\tau)$, $\tau\geq0$,
such that $\omega_{F}(0)=0$ 
and for any $\sfu'_{0},\sfv'_{0}\in\bR$, $x
\in\Omega$, and $\sfu''\in\bS$ we have  
$$
\big|F\big(\sfu'_{0},  \sfu'' ,x\big)-F\big(\sfv'_{0},  \sfu'' ,x\big)\big|\leq
 \omega_{F}\big(|\sfu'_{0}-\sfv'_{0}|\big)|\sfu''|.
$$ 

\end{assumption}

\begin{remark}
                                               \label{remark 1,24,2}
It is useful to note that Assumptions \ref{assumption 10.5.1}
and \ref{assump1} (iv) imply that $F(\sfu'_{0}, 0 ,x)=0$
for any $\sfu'_{0} \in\bR$ and $x
\in\Omega$.
Also
observe that, apart from (iv), Lipschitz continuity in $\sfu''$,
and measurability,
nothing is imposed on $F$ if $|\sfu''|\leq\tau_{0}$.
\end{remark}

 \begin{definition}
                                          \label{definition 1,13,1}
For a function $u\in C(\bar\Omega)$ 
set
$$
\omega_{u}(\Omega,\rho)=
\sup \big\{
 \big|u(x_{1})-u(x_{2})\big|  :x_{1},x_{2}\in\Omega,  
|x_{1}-x_{2}|\leq\rho\big\},
$$
$$
\omega_{F,u,\Omega}(\rho)=\omega_{F}\big(\omega_{u}(\Omega,\rho)\big),
$$
and in the formulations  of
\index{$F$@Functions!$\omega_{u}(\Omega,\rho)$}%
\index{$F$@Functions!$\omega_{F,u,\Omega}(\rho)$}%
 a theorem, lemma,...
  let us say that a certain constant depends only on A,B,..., and
 the function
$\omega_{F,u,\Omega}$ if it depends only on A,B,..., and on the maximal
 solution of an inequality like $N_{0}\omega_{F,u,\Omega}(\rho)\leq
1/2$, where the range of $\rho$ and
 the value of $N_{0}$ depending only on 
A,B,... could be always traced down in our arguments.
 
\end{definition}

To finish the setting, take $d_{0}=d_{0}(d,\delta,\|b\|)
\in(d/2,d)$
from \cite{Kr_20} and take $p\in (d_{0},d)$.
In the statement of the following theorem
we use the function $\bar R(p)$, which is introduced
before Lemma \ref{lemma 5.30.1}
(see \eqref{1.5.1}).

 \begin{theorem}
                                         \label{theorem 11.16.1}
Under the above assumptions
 there exist    constants $\hat\theta,\theta\in (0,1]$, depending
only on $d$, $p$, $\delta$, and  $K_{F}$,
such that, if Assumptions \ref{assump1} and \ref{assumption 10.5.1}
are satisfied with these
 $\theta$ and $\hat\theta$, respectively, then, for any
 $u\in W^{  2}_{ p ,\loc}(\Omega )\cap C(\bar\Omega )$
 that satisfies \eqref{7.29.1} in $\Omega $  (a.e.)
and    $0< \rho<\rho_{\inside}(\Omega)\wedge1\wedge \bar R(p)$, 
where $\rho_{\inside}(\Omega)$
is the interior radius of $\Omega$,
 we have    
\begin{equation}
                                               \label{1,25,3}
\|u\|_{W^{ 2}_p(\Omega  ^{\rho})}\le N
\| \bar G\|_{L_p(\Omega )}+N\rho^{-2}\|u\|_{C( \Omega ) }
 +N \tau_{0} ,
\end{equation}
where the constants $N$  depend only on
 $K_{0}$, $K_{F}$, $d$, $p$, $\delta$, $\|b\|$,  
$R_{0}$, 
  $ \diam(\Omega)$, and the 
function $\omega_{F,u,\Omega}$.      
\end{theorem}

This theorem is proved in Section \ref{section 10.10.2}  
after we develop necessary results in Section \ref{section 10.11.1}.

To state an existence result we
need the following additional assumptions.

\begin{assumption}
                                 \label{assumption 3.11.2}    
 The function   $H(\sfu, x)$ is continuous in $\sfu$
for any $x$, is
  Lipschitz continuous with respect to $\sfu''$, and
$D_{\sfu''}H\in \bS_{\delta}$
at 
all points of differentiability
of $H$ with respect to $\sfu''$.
\end{assumption}

\begin{assumption}
                                 \label{assumption 1.6.1} 
There exists  $n_{0}\geq 0$ such that
for any $x\in\{\bar G>n_{0}\}$ we have 
$D_{\sfu''}F\big(\sfu'_{0},\sfu'',x\big)\in \bS_{\delta}$
at 
all points of differentiability
of $F\big(\sfu'_{0},\sfu'',x\big)$ with respect to $\sfu''$.

\end{assumption}

\begin{assumption}
                                         \label{assumption 3,7,3}
For all values of the  arguments,
\begin{equation}
                                                     \label{3,7,4}
 H(\sfu',0,x)\sign \sfu'_{0}\leq  
\gb(x)\big|[\sfu']\big|+\bar{G}(x)\quad (\sign 0:=\pm 1) .
 \end{equation}
 
\end{assumption}

Here is   our result concerning the solvability
of \eqref{7.29.1} in Sobolev spaces. 
 We fix $p\in(d_{0},d)$ and a function 
$$
g\in C(\partial\Omega).
$$

\begin{theorem}
                                    \label{theorem 10.5.1} 
 There exist    constants $\hat\theta,\theta\in (0,1]$, depending    
only on $d$, $p$, $\delta$, $\|b\|$, and  $K_{F}$,  which
are, generally, smaller than
$\hat\theta,\theta$ from Theorem \ref{theorem 11.16.1} and
such that, if Assumptions \ref{assump1} and \ref{assumption 10.5.1}
are satisfied with these
 $\theta$ and $\hat\theta$, respectively, and
  Assumptions    
 \ref{assumption 3.11.2}, \ref{assumption 1.6.1},   and  
 \ref{assumption 3,7,3} are also satisfied 
and $\bar G\in L_{p}(\Omega)$,
then
  there exists  
$u\in W^{ 2}_{p,\loc}(\Omega)\cap C(\bar\Omega)$
 satisfying \eqref{7.29.1} in $\Omega$  (a.e.)  and such  
that $u=g $ on $\partial\Omega$. Furthermore, in $\Omega$
\begin{equation}
                                                       \label{1.8.2}
|u|   
\leq N\|\bar G\|_{L_{p}(\Omega)}+
\sup_{  \partial \Omega}|g|,
\end{equation}
where $N$ depends only on $p,d,\delta, \|b\|$,
and the diameter of $D$.
\end{theorem}

The proof of this theorem is given in Section \ref{section 1.9.1}.

\begin{remark}
                                            \label{remark 2,16,2}
Since none of characteristics of $\Omega$,
apart from $\rho_{\inside}(\Omega)$ and
  $ \diam(\Omega)$ enters Theorem \ref{theorem 11.16.1},
one can use Theorem \ref{theorem 10.5.1} to prove
the solvability in much worse domains than those
satisfying the exterior ball condition.
Usually one does it by approximating from inside
a given domain, say with smooth ones. For instance,
it would suffice to have 
\begin{equation}
                                                         \label{2,19,4}
\lim_{\rho\downarrow0}\inf_{x\in\partial \Omega}
\frac{\big|B_{\rho}(x)\cap\Omega^{c}\big|}{\rho^{d}}>0,
\end{equation}
see, for instance, Theorem 3.1 of \cite{Sa_88}
or Theorem \ref{theorem 1.6.1}.

We are not pursuing this path and leave  it to the interested
reader.
\end{remark}

\begin{remark}
                                      \label{remark 10.30.1}   
Observe that
generally there is no  uniqueness  
in Theorem \ref{theorem 10.5.1}.
For instance, in the one-dimensional case
the (quasilinear) equation 
$$
 D^{2}u +\sqrt{12|Du |} =0
$$
 for $x\in(-1,1)$
with zero boundary data has two solutions:
one is identically equal to zero and the other one is
 $1-|x|^{3}$.

Another example is given by the (semilinear) equation 
$$
D^{2}u+2u(1+\sin^{2}x+u^{2})^{-1}=0
$$
on $(-\pi/2,\pi/2)$ with zero boundary condition. Again there
are two solutions:
one is $\cos x$ and the other one is identically equal to zero.

To have uniqueness we need different assumptions
(see, for instance, Section 4.1:2 in \cite{Kr_18}).
\end{remark}

\begin{example}
                                 \label{example 5.16.1}
Let $d=3$, $f,\bar G\in L_{p}(\Omega)$,  $\gb\in L_{d}(\Omega)$,
$\alpha\in(0,1]$. Let $w(t)$, $t\in[0,\infty)$,
be a  continuously differentiable
 function  
with sufficiently small derivative.
Then  the   
equation
$$
H(Du,D^{2}u,x):= \bar G(x)\wedge|D_{12}u|
+\bar G(x)\wedge|D_{23}u|+\bar G(x)\wedge|D_{31}u|
 $$
\begin{equation}
                                                      \label{3.4.1}
+ \Delta u+w(|D^{2}u|)+\gb(x)|Du|^{\alpha}-f(x)=0
\end{equation}
satisfies our assumptions and
Theorem \ref{theorem 10.5.1} is applicable.

 Observe that $H$ in \eqref{3.4.1}
is neither convex nor concave with respect to $D^{2}u$.
Also note that we can replace $\Delta u$ with $a^{ij}(x)D_{ij}u$
if $a(x)=(a^{ij}(x))$ is an $\bS_{\delta}$-valued
VMO-function such that $a(x)\geq(\delta^{ij})$.
\end{example}
\begin{example}
                                             \label{example 1,13,1}
Let $A$ and $B$ be some countable sets and
assume that for $\alpha\in A$, $\beta\in B$,   $x\in\bR^{d}$,
  and $\sfu'\in\bR^{d+1}$
we are given an $\bS_{\delta}$-valued function 
 $a^{\alpha }(\sfu'_{0},x)$ (independent of $\beta$)
and a real-valued function $b^{\alpha\beta}(\sfu',x)$.
 Assume that 
  these functions are measurable
in $x$, $a^{\alpha }$ and $b^{\alpha\beta}$ are
 continuous with respect to $\sfu' $  
uniformly with respect to $\alpha,\beta,x$,
and  
$$
\big|b^{\alpha\beta}(\sfu',x)\big|\leq \gb(x)
\big|\big(\sfu'_{1},...,\sfu'_{d}\big)\big|+\bar G(x),
$$
where $\bar G \in L_{p}(\Omega)$ and $\gb\in L_{d}(\Omega)$.
Next assume that there is an $R_{0}\in(0,\infty)$ such that
for any   $z\in\Omega$, $r\in(0,R_{0}]$, and $\sfu'_{0}\in\bR$ one 
can find
$\bar{a}^{\alpha}\in\bS_{\delta}$ (independent of $x$) such that  
$$
 \dashint_{B_{r}(z)}\sup_{\alpha\in A} 
\big|a^{\alpha}\big(\sfu'_{0},x\big)-\bar{a}^{\alpha}\big|\,dx\leq\theta,
$$
where $\theta$ is sufficiently small
(to accommodate Theorem \ref{theorem 10.5.1}).

Consider equation \eqref{7.29.1},
where
$$
H(\sfu,x):=
\infsup_{\beta\in B\,\,\alpha\in A}
\Big[\sum_{i,j=1}^{d}
a^{\alpha }_{ij}\big(\sfu'_{0},x\big)\sfu''_{ij}+
b^{\alpha\beta}(\sfu',x)\Big].
$$
As in Example 10.1.24 of \cite{Kr_18} one easily sees that
Theorem \ref{theorem 10.5.1}  is applicable.
\end{example}

\begin{example}
                                              \label{example 1.21.1}
A further specification of Example \ref{example 1,13,1}
is given by linear equations. Suppose that we are given an $\bS_{\delta}$-valued  measurable function
 $a ( x)$  
and an $\bR^{d}$-valued function $b(x)$ such that $b\in L_{d}(\Omega)$.
 
Next assume that there is an $R_{0}\in(0,\infty)$ such that
for any   ball $B\subset \bR^{d}$ of radius smaller than $R_{0}$
$$
 \dashint_{B }  
 |a(x  )-\bar a_{B} |\,dx\leq\theta,\quad
\bar a_{B}=\dashint_{B }a(x)\,dx.
$$
By using $d,\delta$, and $\|b\|_{L_{d}(\Omega)}$ find
$d_{0}$ as before Theorem \ref{theorem 11.16.1}
and take $p\in(d_{0},d)$. Suppose that we are given
$f\in L_{p}(\Omega)$, nonnegative bounded $c$
on $\Omega$, and $g\in C(\partial\Omega)$.
Consider the equation
$$
a^{ij}D_{ij}u+b^{i}D_{i}u-cu+f=0
$$
in $\Omega$ with boundary condition $u=g$ on 
$\partial\Omega$.

In this situation one can obviously take $F(\sfu'',x)=a^{ij}(x)\sfu''_{ij}$
and satisfy Assumption \ref{assump1}  with $\bar F(\sfu)=
\bar a^{ij}_{B_{r}(z)} \sfu''_{ij}$ and $\tau_{0}=0$.
Assumptions \ref{assumption 10.5.1} (with $\hat \theta=0$,
$K_{0}=\sup c$,
$\gb=|b|$, $\bar G=|f|$), \ref{assumption 3.11.2}, 
\ref{assumption 1.6.1}, and \ref{assumption 3,7,3}
are also satisfied. Therefore, by Theorem \ref{theorem 10.5.1},
if $\theta$ is sufficiently small, depending only on
$d,p,\delta$, and $\|b\|_{L_{d}(\Omega)}$, 
the above boundary value problem has a solution
in $u\in W^{2}_{p,\loc}(\Omega)\cap C(\bar \Omega)$.
Owing to Theorem \ref{theorem 1.14.1} this solution is
unique and this theorem in combination with Theorem 
\ref{theorem 11.16.1}  shows that
for all sufficiently small $\rho>0$
$$
\|u\|_{W^{ 2}_p(\Omega  ^{\rho})}\le N\rho^{-2}\big(
\| f\|_{L_p(\Omega )}+ \|g\|_{C( \partial\Omega ) }\big).
$$
Just in case, observe that how small $\rho$ is depends on the {\em
function\/} $|b|$ and not only on its $L_{d}$-norm.
The main novelty in this example is that $b\in L_{d}(\Omega)$,
and even if $a$ is continuous the result was not known before.

\end{example}

We finish the section with a general comment.
In the proofs of various results  we use
the symbol $N$  to denote finite 
nonnegative constants
which may change from one occurrence to another and
we do not always specify on which data these  constants
depend. In these cases the reader should remember
that, if in the statement of a result there are constants
called $N$ which are claimed to depend only on certain
parameters, then in the proof of the result
the constants $N$ also depend only on the same
parameters unless specifically stated otherwise.
Of course, if we write 
$
N=N(...),
$
 this means that $N$ depends only
on what is inside the parentheses. 
Another point is that when we say that certain constants
depend only on such and such parameters we mean, in particular,
  that the dependence is such that these constants stay bounded
as the parameters vary in compact subsets of their ranges.

 \mysection{Some results from  
\protect\cite{Kr_19_1}  and \protect\cite{Kr_20}}
                                     \label{section 10.10.2}  
The proofs of Theorems \ref{theorem 11.16.1}
 and \ref{theorem 10.5.1}
is based on 
some results from  
\protect\cite{Kr_19_1}  and \protect\cite{Kr_20}
which we collect here.

Let $F(\sfu'')$ be a  {\em convex\/}
function
defined for $\sfu''\in\bS$ such that at all points of 
its differentiability we have
$$
D_{\sfu''}F(\sfu'') \in\bS_{\delta},
 $$
where $\delta\in(0,1]$ is a fixed number. Introduce
$\cL(\delta,\|b\|)$ as the set of operators
$$
L=a^{ij}D_{ij}+b^{i}D_{i},
$$
where $a=(a^{ij})$ is a measurable $\bS_{\delta}$-valued 
function on $\bR^{d}$, $b=(b^{i})$ is a 
measurable $\bR^{d}$-valued function such that
$$
\| b\|_{L_{d}(\bR^{d})}\leq \|b\|.
$$
 
We need the following which for bounded $b$
is found in \cite{Ca_95} and for
  $b\in L_{d+\varepsilon}(\Omega)$
 in \cite{Fo_98}. This is
Corollary 3.1 of \cite{Kr_19_1}. 

\begin{theorem}
                      \label{theorem 1.14.1}
There is a constant $d_{0}=d_{0}(d,\delta,\|b\|)\in(d/2,d)$
such that if
  $p\in[d_{0},\infty)$, $\Omega$ is a bounded domain in $\bR^{d}$, and
$u\in W^{2}_{p,\loc}(\Omega)\cap C(\bar \Omega)$, then for any
  nonnegative measurable function $c$   on $\Omega$
and   $L\in \cL(\delta,\|b\|)$ we have  in $\Omega$
\begin{equation}
                                                 \label{10.8.1}
u\leq N\|(Lu-cu)_{-}\|_{L_{p}(\Omega)}+
\sup_{\partial \Omega}u_{+},
\end{equation}
where $N$ depends only on $p,d,\delta, \| b\|$,
and the diameter of $\Omega$.

\end{theorem}

Here is Theorem 3.2 of \cite{Kr_19_1}, which is useful
while passing to the limit in our nonlinear equations.

\begin{theorem}
                                   \label{theorem 10.14.2}
Let $p\geq d_{0}$, $R\in(0,\infty]$, and 
$L\in\cL(\delta,\|b\|)$. Then there exists
a constant
$N=N(p,d,\delta,\|b\|)\geq0$ such that for any $\lambda>0$    
and $u\in W^{2}_{p,\loc}(B_{R})\cap C(\bar B_{R})$ 
($B_{\infty}=\bR^{d}$, $C(\bR^{d})$
is the set of bounded continuous functions on $\bR^{d}$)
 we have
$$
\lambda\|u_{+}\|_{L_{p}(B_{R/2})}
\leq N\|(\lambda u-Lu)_{+}\|_{L_{p}(B_{R})}
+N\lambda R^{d/p}e^{-R\sqrt{\lambda}/N}\sup_{\partial B_{R}}u_{+},
$$
where the last term should be dropped if $R=\infty$.
\end{theorem}

 We also need the following   Theorem 4.5 of \cite{Kr_19_1},
which is similar to the Fanghua Lin theorem and is used
as one of the main tools in the way the theory
of fully  nonlinear elliptic equations is developed in
\cite{Kr_18}.
\begin{theorem}
                       \label{theorem 10.6.10}
Let  $R\in(0,\infty)$,  $p\in[d_{0},\infty)$, 
$u\in W^{2}_{p,\loc}(B_{R})\cap C(\bar B_{R})$,
$L\in \cL(\delta,\|b\| )$, and
$c\in L_{d_{0}}(B_{R})$, $c\geq0$. Then
$$
\Big(\dashint_{B_{R}}|D^{2}u|^{\gamma } 
\,dx\Big)^{1/\gamma} \leq N \Big(\dashint_{B_{R}}
|L u -cu|^{p}\,dx \Big)^{1/p}
 +NR^{ -2}\sup_{\partial B_{R}}|u|,
$$
where $\gamma=\gamma(  d,\delta,\|b\| )\in(0, 1)$ 
and $N$ depends only on
$d,\delta,\|b\|,p$, and $R^{2-d/d_{0}}\|c\|_{L_{d_{0}}(B_{R} )}$.
\end{theorem}

The following is Corollary 4.11 of \cite{Kr_19_1}
about the boundary behavior of solutions
of linear equations which easily carries over to
the nonlinear case.

\begin{theorem}
                                       \label{theorem 1.6.1}
Let $D$ be a bounded domain in $\bR^{d}$, $0\in\partial D$,
and assume that for some constants $\rho,\gamma>0$ and 
any $r\in (0,\rho)$ we have $|B_{r}\cap D^{c}|\geq \gamma
|B_{r}|$. Suppose that we are given a function $u\in W^{2}_{d_{0},\loc}
(D)\cap  C(\bar D)$ and let  $w(r)$ be a concave 
continuous function on $[0,\infty)$
such that $w(0)=0$ and $|u(x)-u(0)|\leq
w(|x|)$ for all $x\in \partial D$.
Then for $x\in D$ we have
$$
|u(x)-u(0)|\leq N|x|^{\beta}\|Lu\|_{L_{d_{0}}(D)}
+\omega\big(N|x|^{\beta/2}),
$$
where  $L\in \cL(\delta,\|b\| )$ and $N$  depends only on 
$ d,\delta,\|b\|,\gamma,\rho$, and the diameter
of $D$.
\end{theorem}

The following is Corollary 6.8 of \cite{Kr_20}
about estimates of the H\"older constant of solutions.

\begin{theorem}
                                      \label{theorem 1.6.4}
Let $R\in(0,\infty)$, $p\geq d_{0}$,  
and let $u\in W^{ 2}_{p}(B_{2R})$ and $L\in \cL(\delta,\|b\| )$. Define
$f= Lu$. Then there exists a constant $N$, which depends
only on $p,d,\|b\|$, and $\delta$, such that    
$$
\big|u(x_{1})-u(x_{2})\big| \le NR^{-\alpha}
 |x_{1}-x_{2}|^{\alpha}
\big(\sup_{\bar B_{ 2R}}|u|
+R^{2-d/p}\Vert f
\Vert _{L_{p}(B_{ 2R})}\big)
$$
  for $x_{1}$, $x_{2}\in  B_{R}$
with  $\alpha=\alpha(d,\delta,\|b\|)\in(0,1)$.
\end{theorem}

We also need the following result by Safonov
(see \cite{Sa_84}, \cite{Sa_88}, or Section 10.3 in \cite{Kr_18}).
 This is another building block
in the way the theory
of fully  nonlinear elliptic equations is developed in
\cite{Kr_18}.
\begin{theorem}
                                        \label{theorem 10.12.2}
There exists a constant $\alpha_{0}=\alpha_{0}(\delta,d)\in(0,1)$
such that
for any  
 $g\in C(\partial B_{2})$ there exists a unique
$v\in C(\bar{B}_{2})\cap 
C^{2+\alpha_{0}}_{\loc}(B_{2})$ satisfying   
\begin{equation}
                                              \label{10.20.5}
F(D^{2}v)=0\quad\text{in}\quad B_{2},\quad v=g\quad
\text{on}\quad\partial B_{2}.
\end{equation}
Furthermore,
$$
\big|D^{2}v(x)-D^{2}v(y)\big|\leq N|x-y|^{\alpha_{0}}
\sup_{\partial B_{2}}|g-p|
$$
as long as $x,y\in B_{1}$,
where $p$ is an arbitrary polynomial
of degree 2 on $\bR^{d}$ and $N$ depends only on $\delta$ and $d$.
\end{theorem}

Below   we fix
$
\alpha\in(0, \alpha_{0}]$.
Here is Lemma 10.3.2 of \cite{Kr_18}.

\begin{lemma}
                                          \label{lemma 12.16.01}

Let $r\in(0,\infty)$, $\nu\geq2$ and let $\phi\in
C(\partial B_{\nu r})$. Then there exists
a unique $v\in C(\bar{B}_{\nu r})\cap C^{2+\alpha}_{\loc}
(B_{\nu r})$ such that
$$
F(D^{2}v)=0\quad\text{in}\quad B_{\nu r},\quad v=\phi\quad
\text{on}\quad\partial B_{\nu r}.
$$
Furthermore,   
$$
\dashint_{B_{r}}\dashint_{B_{r}}\big|D^{2}v(x)-D^{2}v(y)\big|
\,dxdy\leq N(d,\alpha,\delta)\nu^{-2-\alpha} r^{-2}\sup_{\partial
B_{\nu r}}|\phi|.
$$
\end{lemma}

Finally, we will use the following, which allows us to use
a version of the Fefferman-Stein theorem.
 
\begin{lemma}
                                      \label{lemma 10.5.2}   
 
Let $r\in(0,\infty)$ and $\nu\in[2,\infty)$. Then for any
$u\in W^{2}_{d_{0} }(B_{\nu r})$    we have   
$$
\Big(\dashint_{B_{r} }\dashint_{B_{r} }
\big|D^{2}u(x)-D^{2}u(z)\big|^{\gamma}\,dxdy\Big)^{1/\gamma}
$$
\begin{equation}
                                           \label{10.6.1}
\leq N\nu^{d/\gamma}
\Big(\dashint_{B_{\nu r} }
\big|F[u]\big|^{d_{0}}\,dx\Big)^{1/d_{0}}+
N\nu^{-\alpha}\Big(\dashint_{B_{\nu r} }|D^{2}u|^{d_{0}}\,dx
\Big)^{1/d_{0}},
\end{equation}
where $N$ depends only on $d $, $\delta$,   and $\|b\|$
and $\gamma$ is taken from Theorem \ref{theorem 10.6.10}.
\end{lemma}

Proof. Define
$v$ to be a unique $C(\bar{B}_{\nu r} )
\cap C^{2+\alpha}_{\loc}(B_{\nu r} )$-solution of
the equation $F[v]=0$ in $B_{\nu r}$ with boundary condition
$v=u$ on $\partial B_{\nu r} $. Such a function exists
by Lemma \ref{lemma 12.16.01}.
Furthermore, $v(x)-b^{i} x^{i} -c$ satisfies
 the same equation
for any constants $b^{i},c$. Hence by Lemma \ref{lemma 12.16.01}
and H\"older's inequality
$$
I_{r} :=\Big(
\dashint_{B_{r} }\dashint_{B_{r} }\big|D^{2}v(x)-D^{2}v(y)\big|^{\gamma}
\,dxdy\Big)^{1/\gamma}
$$
$$
\leq N \nu^{-2-\alpha} r^{-2}\sup_{x\in\partial
B_{\nu r} }\big|u(x)-(D_{i}u)_{B_{\nu r}} x^{i} 
-u_{B_{\nu r}}\big|.
$$
By Poincar\'e's inequality (recall that $d_{0}>d/2$)
 the last supremum is dominated
by a constant times
$$
\nu^{2} r^{2}\Big(\dashint_{B_{\nu r} }|D^{2}u|^{d_{0}}\,dx 
\Big)^{1/d_{0}}.
$$
It follows that
\begin{equation}
                                            \label{10.6.3}
I_{r} \leq N\nu^{-\alpha}\Big(\dashint_{B_{\nu r} }|D^{2}u|^{d_{0}}\,dx
\Big)^{1/d_{0}}.
\end{equation}

Next, the function $w=u-v$ is of class $C(\bar{B}_{\nu r} )
\cap W^{2}_{d_{0},\loc}(B_{\nu r} )$ and for an operator  
 $ L\in\cL_{\delta,0 }$ we have 
$$
F[u]-F[v]= L(u-v),\quad
   L w=F[u] 
$$
in $B_{\nu r}$ (a.e.). Moreover, $w=0$ on $\partial B_{\nu r}$.
Therefore,
by Theorem \ref{theorem 10.6.10}
$$
\dashint_{B_{r} }|D^{2}w|^{\gamma}\,dx
\leq \nu^{d}\dashint_{B_{\nu r} }|D^{2}w|^{\gamma}\,dx
\leq N\nu^{d}\Big(\dashint_{B_{\nu r} }
\big|F[u]\big|^{d_{0}}\,dx\Big)^{\gamma/d_{0}}.
 $$

Upon combining this result with \eqref{10.6.3} we come 
to \eqref{10.6.1} and the lemma is proved.   \qed

\mysection{ Proof of Theorem \protect\ref{theorem 11.16.1}}
                                 \label{section 10.11.1}

Here we suppose that Assumptions  
\ref{assumption 10.5.1}  and \ref{assump1} are satisfied
with $\theta$, $\hat \theta$ to be specified later.
Thus, we suppose that all assumptions stated before Theorem
\ref{theorem 11.16.1}
are satisfied. 

First we recall the following Lemma 10.4.1 of \cite{Kr_18}.

\begin{lemma}
                                         \label{lemma 10.6.4}
 For any $q\in[1,\infty)$ and $\mu>0$
there is a $\theta=\theta(d,\delta,K_{F},\mu,q)>0$ such that,
if Assumption \ref{assump1} is satisfied with this $\theta$,
then the following holds:

 for any $\sfu'_{0}\in\bR$, $r\in(0,R_{0}]$ and $z\in\Omega$ such that
$B_{r}(z)\subset\Omega$ we have   
$$
 \dashint_{  B_{r}(z)}\sup_{\substack{\sfu''\in\bS,\\|\sfu''|>
\tau_{0}}}
\frac{\big|F\big(\sfu'_{0},\sfu'',x\big)-\bar{F}(\sfu'')\big|^{q}}{|\sfu''|^{q}}\,dx\leq \mu^{q},
$$
where $\bar F=\bar F_{z,r,\sfu'_{0}}$.  
\end{lemma}

Below $\gamma$ is taken from Theorem \ref{theorem 10.6.10}.

\begin{lemma}
                                       \label{lemma 19.12.1}
Let $r\in(0,\infty)$ and
 $\nu\geq 2$ be such that  $\nu r\leq 
 R_{0} $ and $\Omega^{\nu r}\ne\emptyset$.
Take   
$$
\mu\in(0,\infty),\quad \beta\in(1,\infty),
$$ 
and suppose that the assertion of Lemma \ref{lemma 10.6.4} holds
with $q=\beta d_{0}$.
Take   a function $u\in W^{ 2}_{d_{0} }(\Omega )$,
and 
for $x_{0} \in\Omega^{\nu r}$  denote  
$$
I_{r}(x_{0})=\Big(\dashint_{B_{r}(x_{0}) }
\dashint_{B _{r}(x_{0}) }
\big|D^{2}u(x_{1})-D^{2}u(x_{2})\big|^{\gamma}\,dx_{1}dx_{2}\Big)^{1/\gamma}.
$$

Then     for any $x_{0}\in\Omega^{\nu r}$  
\begin{align}
I_{r} (x_{0})\leq &\, N\nu^{d /\gamma}\Big(
\dashint_{B_{\nu r }(x_{0}) } \big|F[u]\big|^{d_{0}}
 \,d x\Big)^{1/d_{0}}+N\tau_{0}\nu^{d /\gamma}
\nonumber \\
                                             \label{10.9.01}
&\,+N\Big[\big(\mu+\omega_{F,u,\Omega}(\nu r )\big)
\nu^{d /\gamma} +\nu^{-\alpha}
\Big]
\Big(
\dashint_{B_{\nu r }(x_{0}) }
|D^{2}u |^{\beta' d_{0}}\,dx\Big)^{1/(\beta'd_{0})},
\end{align}
where $\beta'=\beta/(\beta-1)$  and 
$N$ depends only on $ d,K_{F}$, 
$\delta$, and $\|b\|$.

\end{lemma}

This lemma is proved in the same way as Lemma 10.4.2 of \cite{Kr_18},
basically, using only H\"older's inequality and
Lemmas \ref{lemma 10.6.4} and \ref{lemma 10.5.2}.
By the way the term $\omega_{F,u,\Omega}(\nu r )$
appears because of Assumption \ref{assump1} (iv).

Lemma \ref{lemma 19.12.1}
allows us to follow the proof
of Lemma 10.4.3 of \cite{Kr_18}, which we prefer here
to split into two parts. Here is the first part.

\begin{lemma}
                                       \label{lemma 1.4.1}    
Take $p\in(d_{0},d)$,
   $R\in(0,1]$,  and  $u\in W^{2}_{p}(B_{2R })$.
Take $\mu\in(0,\infty)$ and suppose that the assertion 
of Lemma \ref{lemma 10.6.4} holds with $q=\beta p$,
where $\beta$ is so large that $\beta' d_{0}<p$. Take $\varepsilon\in(0,1]$
  and
 let $0<R_{1}<R_{2}\leq 2R$ be such that 
\begin{equation}
                                                \label{10.14.1}
 R_{2}-R_{1}\leq \varepsilon R_{0},\quad R_{2}\leq  2R_{1}.
 \end{equation}
Assume that $B_{2R }\subset\Omega$.
 Then
there exist    constants  $N$, $N_{1}$, and  $N_{2}$,
depending only on   $d$, $p$,   $K_{F}$,  $\delta$, 
$\beta$, and  $\|b\|$,
  such that  
\begin{align}
\|D^{2}u\|_{L_{p}(B_{R_{1}})}\leq &\,
N_{1} \big\|F[u]\big\|_{L_{p}(B_{R_{2} })}
+N\tau_{0}  R_{1} ^{d /p}
\nonumber\\ 
&\, +\Big[ N_{2}\big(\mu+\omega_{F,u,B_{2R}}
(\varepsilon R_{0})\big)+1/16\Big] 
\|D^{2}u\|_{L_{p}(B_{R_{2} })} 
\nonumber\\ 
                                                      \label{4,21,1}
&\, +N (R_{2}-R_{1})^{-\chi_{1}}R_{1}^{-\chi_{2}+\chi_{1}}
\big\|\,|D^{2}u|^{\gamma}
\big\|^{1/\gamma}_{L_{1}(B_{2R })}, 
\end{align}
where 
\begin{equation}
                                                    \label{1.4.3}
\chi_{1}=(d +2)/\gamma,\quad \chi_{2}=d /\gamma-d /p.
\end{equation}

\end{lemma}

Proof. Proof. For   $\rho>0$  and $x\in
\bR^{d}$ introduce 
$$
h^{\sharp}_{ \gamma,\rho}(x) =
\sup_{\substack{r\in(0,\rho],\\
B_{r}(x_{0})\ni x}}\Big(\dashint_{B_{r}(x_{0}) }
\dashint_{B_{r}(x_{0}) }
\big|h(x_{1})-h(x_{2})\big|^{\gamma}\,dx_{1}dx_{2}\Big)^{1/\gamma},
$$
\begin{equation}
                                                    \label{7,2,1}
\bM  h(x) =\sup_{\substack{r>0,\\
B_{r}(x_{0})\ni x}}\dashint_{B_{r}(x_{0}) }
|h(y)|\,dy,
\end{equation}
\index{$F$@Functions!$h^{\#}_{ \gamma,\rho}$}%
\index{$C$@Operators!$\bM  h$}%
whenever these definitions make sense.
 
Then
take $\nu\geq2$ and set      
$$
 r_{0}=(R _{2}-R_{1})/(\nu+1) .
$$

Next, take $x,x_{0}$, and $r>0$ such that 
$$
r\leq r_{0},\quad
 x\in B_{R_{1}},\quad x\in B_{r }(x_{0})
$$
and observe that, since $R_{2} -\nu r_{0}=R_{1}
+r_{0}$, we have $x_{0}\in B_{R_{2}-\nu r_{0}}$
and $B_{\nu r}(x_{0})\subset B_{R_{2}}$. Also $\nu r\leq
\nu r_{0}\leq   R_{0}$.
Therefore, by Lemma \ref{lemma 19.12.1} applied to $\Omega=B_{R_{2} }$,
    we have  
 (note $x_{0}$ on the left and $x$ on the right) 
\begin{align*}
I_{r}(x_{0})\leq &\, N\nu^{d /\gamma}\bM^{1/d_{0}}
\big(\big|F[u]\big|^{d_{0}}I_{
B_{R_{2} }}\big)(x)+N\tau_{0}\nu^{d /\gamma}
\\ 
&\, +N\Big[\big(\mu+\omega_{F,u,B_{2R}}(\nu r_{0})\big)
\nu^{d /\gamma} +\nu^{-\alpha}\Big]
\bM^{1/(\beta'd_{0})}\big(|D^{2}u |^{\beta'd_{0}}I_{
B_{R_{2} }}\big)(x)
\end{align*}
with $N$ depending only on   $ d,K_{F}$, 
and $\delta$.
It follows that in $B_{R_{1}}$  
\begin{align*}
(D^{2}u)^{\sharp}_{\gamma,r_{0}} \leq &\, N\nu^{d /\gamma}
\bM^{1/d_{0}}\big(\big|F[u]\big|^{d_{0}}I_{
B_{R_{2} }}\big) +N\tau_{0}\nu^{d /\gamma}
\\ 
&\,+N\Big[\big(\mu+\omega_{F,u,B_{2R}}(\varepsilon R_{0})\big)
\nu^{d /\gamma} +\nu^{-\alpha} 
\Big]
\bM^{1/(\beta'd_{0})}\big( |D^{2}u |^{\beta'd_{0}}I_{
B_{R_{2} }}\big) .
\end{align*}

By Theorem C.2.6 of \cite{Kr_18} (which 
is similar to the Fefferman-Stein theorem)
 with   $\kappa=r_{0}/R_{1} \leq
1/3 $ and $\chi_{1},\chi_{2}$ from \eqref{1.4.3}
and the Hardy-Littlewood maximal function theorem
(recall that
  $p>\beta'd_{0}$), we obtain  
\begin{align}
\|D^{2}u\|_{L_{p}(B_{R_{1}})}\leq &\,
N \nu^{d /\gamma}\big\|F[u]\big\|_{L_{p}(B_{R_{2} })}
+N\tau_{0}\nu^{d /\gamma} R_{1} ^{d /p}
\nonumber\\ 
&\, +\Big[N \big(\mu+\omega_{F,u,B_{2R}}
(\varepsilon R_{0})\big)\nu^{d /\gamma}
 +N_{0}\nu^{-\alpha}) \Big]
\|D^{2}u\|_{L_{p}(B_{R_{2} })}
\nonumber\\ 
                                                          \label{6.7.01}
&\, +N\nu^{\chi_{1}}(R_{2}-R_{1})^{-\chi_{1}}R_{1}^{-\chi_{2}+\chi_{1}}
\big\|\,|D^{2}u|^{\gamma}
\big\|^{1/\gamma}_{L_{1}(B_{2R })},  
\end{align}
where   the constants $N$, $N_{i}$ depend only on    
  $ d$, $p$, $K_{F}$, $\|b\|$,
and $\delta$.  
Now we take and fix $\nu\geq2$ so that
$$
N_{0}\nu^{-\alpha}\leq 1/16.
$$
Then \eqref{6.7.01} becomes \eqref{4,21,1}. The lemma is proved. \qed
 
The constant $N_{1}$ in \eqref{4,21,1} depends only
on $d$, $p$, $\beta$,  $K_{F}$, $\|b\|$, and  $\delta$,
and $\beta$ can be easily made to depend only on
$p$ and $d_{0}$. Therefore, the constant $N_{1}$ 
in \eqref{4,21,1} depends only
on $d$, $p$,  $K_{F}$, $\|b\|$, and  $\delta$: $N_{1}=
N_{1}(d,p,K_{F},\|b\|,\delta)$. 
Another constant we need to proceed is the following. For $p\in[1,d)$
and $q=pd/(d-p)$ by interpolation inequalities
there is a constant $N(p,d)$ such that for any $R\in(0,1]$
and $u\in W^{2}_{p}(B_{R})$ we have
$$
\|Du\|_{L_{q}(B_{R})}\leq N(p,d)\|D^{2}u\|_{L_{p}(B_{R})}
+N(p,d)R^{-2}\| u\|_{L_{p}(B_{R})}.
$$
Now, we define $\bar R(p)$ by   requiring
that  $\bar R(p)\in(0,1]$ and 
for any $x\in\Omega$, and
for any $R\in(0,\bar R(p)]$
\begin{equation}
                                              \label{1.5.1}
N_{1}(d,p,K_{F},\|b\|,\delta)\|\gb\|_{L_{d}(B_{2 R}(x))}N(p,d)\leq 1/8.
\end{equation}
\begin{lemma}
                                       \label{lemma 5.30.1}    
Take $p\in(d_{0},d)$,
   $R\in(0,\bar R(p)]$,  and  $u\in W^{2}_{p}(B_{2R })$.
Assume that $B_{2R }\subset\Omega$.
 Then
there exist    constants $\hat\theta,\theta\in (0,1]$, depending  
only on $d$, $p$, $\delta$, $\|b\|$,  and  $K_{F}$,
such that, if Assumptions \ref{assumption 10.5.1} and
 \ref{assump1} are satisfied with these $\hat\theta$ and $\theta$, 
respectively,
  then   
  there is a constant  $N$,
depending only on $R_{0}$, $d$, $p$, $K_{0}$, $K_{F}$,   $\delta$,  
  $\|b\|$, 
and the function $\omega_{F,u,  B_{2R}}$,  such that   
\begin{align}
\|D^{2}u\|_{L_{p}(B_{R })}\leq &\,
N \big\|H[u]\big\|_{L_{p}(B_{2R })}+N\| \bar G\|_{L_{p}(B_{2R })}
+N\tau_{0}    R  ^{d/p}
\nonumber\\[10pt]
                                                      \label{6.16.3}
&\,+N  R^{d/p-d/\gamma }\big\|
\,|D^{2}u|^{\gamma}
\big\|^{1/\gamma}_{L_{1}(B_{2R })}+NR^{-2}\| u\|_{L_{p}(B_{2R })},
\end{align} 
\begin{align}
\|D^{2}u\|_{L_{p}(B_{R})}\leq &\, N\tau_{0} R  ^{d/p}  
+NR^{d/p-2}\sup_{B_{2R}}|u|
\nonumber\\[10pt]
                                             \label{5.30.1}
&\,+N \big( \big\|H[u]\big\|_{L_{p}(B_{2R})}+\|
 \bar G\|_{L_{p}(B_{2R})}\big).
\end{align} 

\end{lemma}

Proof. Take $\varepsilon\in(0,1]$
  and
 let $0<R_{1}<R_{2}\leq 2R$ be as in Lemma \ref{lemma 1.4.1}.
Also take $\mu\in (0,\infty)$ and suppose that
Assumption \ref{assump1} holds with
$\theta=\theta(d,\delta,K_{F},\mu,\beta d_{0})$
(see Lemma \ref{lemma 10.6.4}),
where $\beta=\beta(d_{0},p)$ is so large that $\beta' d_{0}<p$.
Then \eqref{4,21,1} holds. We estimate
$F[u]$ by observing that    
$$
\big|F[u]\big|\leq\big|H[u]\big|+K_{0}|u|+\gb|Du| +\bar G+\hat\theta|D^{2}u|
$$
and that by H\"older's and interpolation inequalities   
with $q=pd/(d-p)$ and by \eqref{1.5.1}
$$
 N_{1} \|\gb D u\|_{L_{p}(B_{R_{2} })}
\leq N_{1}\|\gb \|_{L_{d}(B_{R_{2}})}\|D u\|_{L_{q}(B_{R_{2} })}
$$
$$
\leq (1/8)\|D^{2}u\|_{L_{p}(B_{R_{2} })}
+
(1/8)R_{2}^{-2} 
\|u\|_{L_{p}(B_{R_{2} })}.
$$

Then we take $\hat \theta$ and $\mu$ so small that
$$
N_{1}\hat\theta\leq 1/8,\quad N_{2}\mu\leq 1/8,
$$
and, finally, take 
the largest $\varepsilon\leq1$ such 
that   
$$
N_{2}
\omega_{F,u,B_{2R}}(\varepsilon R_{0})\leq 1/8.
$$
This $\varepsilon$, which depends only on
 $ d$, $p$, $K_{F}$, $R_{0}$, the function $\omega_{F,u,B_{2R}}$,
$\|b\|$, and $\delta$, will appear later
in our arguments and this is the way how the constant $N$ in the statement
of the lemma depends on $\omega_{F,u,B_{2R}}$.

We require   Assumptions \ref{assumption 10.5.1}
and  \ref{assump1}
be satisfied with the above chosen $\hat\theta$ and
 $\theta=\theta(d,\delta,K_{F},\mu,\beta d_{0})$,
 respectively. By combining the above,
 we get 
$$
\|D^{2}u\|_{L_{p}(B_{R_{1}})}\leq
N \big\|H[u]\big\|_{L_{p}(B_{R _{2}})}
+N\tau_{0}   R  ^{d /p}
+(5/8)
\|D^{2}u\|_{L_{p}(B_{R_{2} })}
$$
$$
+N(R_{2}-R_{1})^{-\chi_{1}}R_{1}^{-\chi_{2}+\chi_{1}}\big\|\,|D^{2}u|^{\gamma}
\big\|^{1/\gamma}_{L_{1}(B_{2R })}+NR^{-2}_{2}\|u\|_{L_{p}(B_{2R  })}
+N\| \bar G\|_{L_{p}(B_{2R })}.
 $$
 
  Now we are going to iterate this
estimate by defining $R_{1}=R$ and for $k\geq 1$  
$$
R_{k+1}=R_{k }+c R (n_{0}+k)^{-2},
$$
where the constant $c=O(n_{0})$ is chosen so that $R_{k}
\uparrow 2R$ as $k\to\infty$, that is  
$$
c\sum_{k=1}^{\infty}(n_{0}+k)^{-2}=1,
$$
and 
$n_{0}>0$ is chosen so that for $k\geq 1$  
$$
R_{k+1}-R_{k}=c R (n_{0}+k)^{-2}\leq  R cn^{-2}_{0}
\leq R\leq  R_{k},
$$
which is satisfied if $n_{0}$ is just an appropriate
absolute
constant, and  
$$
R_{k+1}-R_{k}=cR(n_{0}+k)^{-2}\leq  cn^{-2}_{0}\leq \varepsilon R_{0}
$$
(this time we need $n^{-1}_{0}=o(
\varepsilon R_{0} )$ if $\varepsilon  R_{0}\to0$).
Also observe that $R\leq R_{k}\leq 2R$ and 
$$
(R_{k+1}-R_{k})^{-\chi_{1}} R_{k}^{-\chi_{2}+\chi_{1}}
\leq N(n_{0}+k)^{2\chi_{1}} R ^{-\chi_{2}} .
$$

Then for $k\geq1$ we get 
$$
\|D^{2}u\|_{L_{p}(B_{R_{k}})}\leq
N \big\|H[u]\big\|_{L_{p}(B_{R _{k+1}})}
+N\tau_{0}   R  ^{d/p}
+(5/8)
\|D^{2}u\|_{L_{p}(B_{R_{k+1} })}
$$
$$
+N (n_{0}+k)^{2\chi_{1}}R^{-\chi_{2}}\big\|
\,|D^{2}u|^{\gamma}
\big\|^{1/\gamma}_{L_{1}(B_{2R })}+NR^{-2}\| u\|_{L_{p}(B_{2R })}
+N\| \bar G\|_{L_{p}(B_{2R })},
$$
where and below the constants $N$ are as in the statement of 
the lemma.
We multiply both parts of this inequality by $(5/8)^{ k}$
and sum up the results over $k=1,2,...$. Then we cancel
 the  like terms
$$
\sum_{k=2}^{\infty}(5/8)^{ k}\|D^{2}u\|_{L_{p}(B_{R_{k}})},
$$
which are finite since $u\in W^{2}_{p}(B_{2R})$,
and finally take into account that    
$$
 \sum_{k=2}^{\infty}(5/8)^{ k}(n_{0}+k)^{2\chi_{1}}
\leq N n_{0}^{2\chi_{1}}\sum_{k=2}^{\infty}(5/8)^{ k}
+N\sum_{k=2}^{\infty}(5/8)^{ k}k ^{2\chi_{1}}\leq N.
$$
Then we come to \eqref{6.16.3}.

To derive \eqref{5.30.1} observe that, given $u\in W^{2}_{p}(B_{2R})$,
there exists $L\in\cL(\delta,\|b\|)$ such that
$$
H[u]=[H(D^{2}u,Du,u)-H( 0,Du,u)]+H(0,Du,u)=Lu+f,
$$
where $|f|\leq K_{0}|u|+\bar G$. Therefore, $Lu=H[u]-f$
and by Theorem \ref{theorem 10.6.10}
$$
R^{d/p-d/\gamma }\big\|
\,|D^{2}u|^{\gamma}
\big\|^{1/\gamma}_{L_{1}(B_{2R })}\leq
N\|Lu\|_{L_{p}(B_{2R})}+NR^{d/p-2}\sup_{B_{2R}}|u|
$$
$$
\leq N\|H[u],\bar G\|_{L_{p}(B_{2R})}+NR^{d/p-2}\sup_{B_{2R}}|u|,
$$
where we used that $\|u\|_{L_{p}(B_{2R})}\leq NR^{d/p-2}
\sup_{B_{2R}}|u|$
because $R\leq 1$. 
Finally, observing that $\|u\|_{L_{p}(B_{2R})}\leq NR^{d/p}
\sup_{B_{2R}}|u|$ we come from \eqref{6.16.3} to
\eqref{5.30.1} and the lemma is proved.  
\qed

{\bf Proof of Theorem \ref{theorem 11.16.1}}. 
We take the
  constants $\hat\theta,\theta\in (0,1]$ from Lemma \ref{lemma 5.30.1}. 
By that lemma, if $\rho\in(0,\bar R(p)]$
 and $\Omega^{2\rho}\ne\emptyset$ and
$z\in\Omega^{2\rho}$, we have  
$\bar B_{2\rho}(z)\subset\Omega$ and  
$$
\|D^{2}u\|^{p}_{L_{p}\big(B_{\rho}(z)\big)}  \leq
 N\tau_{0}^{p}|\rho |^{d }  +N\rho^{d -2p}\sup_{B_{2\rho}(z)}|u|^{p}
+N\| \bar G\|_{L_{p}\big(B_{2\rho}(z)\big)}^{p}.
 $$
This,  Lemma 10.4.4 of \cite{Kr_18}, implies that,
for $0<3\rho<\rho_{\inside}(\Omega)\wedge3$, 
\begin{equation}
                                                 \label{11.5.1}   
\int_{\Omega^{3\rho}} \big|D^{2}u(x)\big|^{p}\,dx\leq N
\int_{\Omega }  \big|\bar G(x)\big|^{p}\,dx+N\tau_{0}^{p}
 +N \rho^{ -2p}\sup_{\Omega}|u|^{p}.
\end{equation}
 
Using interpolation inequalities also allows us to estimate
the $L_{p}(\Omega^{3\rho})$-norm of $Du$. The theorem is proved. \qed

\mysection{Proof of Theorem \protect\ref{theorem 10.5.1}}

                                                 \label{section 1.9.1}
We give the proof of Theorem \ref{theorem 10.5.1}
after some preparations.
First we use the solvability result from \cite{Kr_18}
in which, however, $\gb$ was assumed to be bounded
and $\bar G\in L_{q}(\Omega)$ for some $q>d$.
Recall that $\bar R(p)$ is introduced in
\eqref{1.5.1}.
\begin{lemma}
                                          \label{lemma 1.6.1}
Assume that $\bar G\in L_{q}(\Omega)$ for some $q>d$
and let $p\in(d_{0},d)$. For $n=1,2,...$ introduce
\begin{equation}
                                                \label{1.8.10}
H^{n}(\sfu,x)=H(\sfu'_{0},n[\sfu']/(n+\gb(x)),\sfu'',x).
\end{equation}
Then there exist constants $\hat\theta,\theta\in (0,1]$, depending
only on $d$, $p$, $\delta$, and  $K_{F}$,
such that, if Assumptions \ref{assump1} and \ref{assumption 10.5.1}
are satisfied with these
 $\theta$ and $\hat\theta$, respectively, then
  for any $n$ there exists a solution  $u_{n}
\in W^{ 2}_{p,\loc}(\Omega)\cap C(\bar\Omega)$
of  the equation
\begin{equation}
                                                \label{1.6.2}
 H^{n}[u_{n}]=0
\end{equation}
(a.e.) in $\Omega$ with boundary data $u_{n}=g$ on 
$\partial\Omega$. Furthermore, for any
  $0< \rho<\rho_{\inside}(\Omega)\wedge1\wedge \bar R(p)$, 
 we have    
\begin{equation}
                                               \label{1,25,30}
\|u_{n}\|_{W^{ 2}_p(\Omega  ^{\rho})}\le N
\| \bar G\|_{L_p(\Omega )}+N\rho^{-2}\|u_{n}\|_{C( \Omega ) }
 +N \tau_{0} ,
\end{equation}
where the constants $N$  depend only on
 $K_{0}$, $K_{F}$, $d$, $p$, $\delta$, $\|b\|$,  
$R_{0}$, 
  $ \diam(\Omega)$, and the 
function $\omega_{F,u,\Omega}$.  
\end{lemma}

Proof. Observe that owing to \eqref{1,14,2}
we have
$$
|G(\sfu'_{0},n[\sfu']/(n+\gb),\sfu'',x)|\leq
\hat\theta|\sfu''|+K_{0}|\sfu'_{0}|+\gb_{n}|[\sfu']|
+\bar G,
$$
and $\gb_{n} = n\gb/(n+\gb)$ is bounded.
Therefore, by Theorem 10.1.14 of \cite{Kr_18}
there exist constants $\hat\theta,\theta\in (0,1]$, depending
only on $d$,   $\delta$, and  $K_{F}$,
such that, if Assumptions \ref{assump1} and \ref{assumption 10.5.1}
are satisfied with these
 $\theta$ and $\hat\theta$, respectively, then
equation \eqref{1.6.2} with given boundary data
has a solution $u_{n} \in W^{ 2}_{q,\loc}(\Omega)\cap C(\bar\Omega)$.
By reducing $\theta$ and $\hat\theta$ in order to accommodate
those in
Theorem \ref{theorem 11.16.1} we prove the second statement.
The lemma is proved. \qed

Now, naturally we want to sent $n\to\infty$.
For a function $u=u(x)$, for which $Du(x)$ is well defined
  we set
$$
H_{u,Du}(\sfu'',x)=H(u(x),Du(x),\sfu'',x).
$$
 
\begin{lemma}
                                          \label{lemma 1.7.1}
Let $p\in(d_{0},\infty)$, $R\in(0,\infty)$,
$u,u_{n}\in W^{2}_{p}(B_{R})$, $n=1,2,...$. Suppose that
$$
M:=\sup_{n}\|u_{n}\|_{W^{2}_{p}(B_{R})}<\infty,
\quad u_{n}\to u\quad\text{weakly in}\quad W^{2}_{p}(B_{R}).
$$
Then there is a subsequence $n'\to\infty$ 
 such that in $B_{R}$
(a.e.)
\begin{equation}
                                                  \label{1.7.2}
\lim_{n'\to\infty}H^{n'}_{u_{n'} ,Du_{n'} }( D^{2}u(x),x)
=H[u](x),
\end{equation}
\begin{equation}
                                                  \label{1.7.30}
\sup_{n'}|H^{n'}_{u_{n'} ,Du_{n'} }( \sfu'',x)|
\leq N(d,\delta)|\sfu''|+\hat H,
\end{equation}
where  
the nonnegative $\hat H$ is such that
$$
\|\hat H\|_{ L_{d_{0}}(B_{R})}\leq
N(d,d_{0},p,R)(K_{0}+\|b\|)(M+1)+\|\bar G\|_{ L_{d_{0}}(B_{R})}.
$$
\end{lemma}

Proof. Let $q=d_{0}d/(d-d_{0})$.  By embedding
theorems $u_{n}\to u$ strongly in $W^{1}_{q}(B_{R})$ and
there exists a subsequence, identified for simplicity
with the original one, such that
$$
\|u_{n+1}-u_{n}\|_{W^{1}_{q}(B_{R})}\leq 2^{-n}.
$$
Then $u_{n },Du_{n }\to
u,Du$ (a.e.) in $B_{R}$ and 
\eqref{1.7.2} follows since $H$ is continuous in $\sfu'$.

Next set
$$
w_{0}=\sum_{n}|u_{n+1}-u_{n}|+|u_{1}|,
\quad w_{1}=\sum_{n}|Du_{n+1}-Du_{n}|+|Du_{1}|.
$$
We have that $w_{0},w_{1}\in L_{q}(B_{R})$, $|u|,|u_{n}|\leq w_{0}$,
$|Du|,|Du_{n}|\leq w_{1}$,
so that
$$
|H^{n}_{u_{n } ,Du_{n } }( \sfu'',x)|\leq N(d,\delta)|\sfu''|+
K_{0}w+\gb w_{1}+\bar G.
$$
This implies \eqref{1.7.30} because by H\"older's
inequality $\gb w_{1}\in L_{d_{0}}(B_{R})$.
The lemma is proved.

To pass to the limit as $n\to\infty$ under the sign of
$H$ which is nonlinear we use the following
replacement of nonlinear operators with linear ones.

\begin{lemma}
                                          \label{lemma 1.7.2}
Let $p\in[1,\infty)$, $u\in W^{2}_{p,\loc}(\Omega)$ satisfy
 \eqref{7.29.1} in $\Omega$  (a.e.). Then there exists
an $\bS_{\delta}$-valued measurable function $a$, $\bR^{d}$-valued
measurable $b$   such that $|b|\leq\gb$ and
(a.e.) in $\Omega$
\begin{equation}
                                                \label{1.7.3}
|a^{ij}D_{ij}u+b^{i}D_{i}u|\leq   K_{0}|u|+\bar G.
\end{equation}
Furthermore, if $u\geq0$ in $\Omega$, then (a.e.) in $\Omega$
(with perhaps different $b$)
\begin{equation}
                                                \label{1.7.4}
 a^{ij}D_{ij}u+b^{i}D_{i}u +\bar G\geq 0
\end{equation}
and  if $u\leq0$ in $\Omega$, then (a.s.) in $\Omega$
\begin{equation}
                                                \label{1.7.5}
 a^{ij}D_{ij}u+b^{i}D_{i}u -\bar G\leq 0.
\end{equation}
\end{lemma}

Proof. By using Assumption \ref{assumption 3.11.2} we get
$$
0=H[u]=H[u]-H(u,Du,0,x)+H(u,Du,0,x)=
a^{ij}D_{ij}u+H(u,Du,0,x),
$$
where by Assumption \ref{assumption 10.5.1}
$$
|H(u,Du,0,x)|\leq K_{0}|u|+\gb|Du|+\bar G=K_{0}|u|+\gb e^{i}D_{i}u
+\bar G,
$$
where $e=Du/|Du|$, which implies that for a function $t(x)$ with values in $[-1,1]$
$$
 H(u,Du,0,x)= t K_{0}|u|+t\gb e^{i}D_{i}u +t\bar G,
$$
and this yields \eqref{1.7.3}.

To prove \eqref{1.7.4} we use the information
provided by Assumption \ref{assumption 3,7,3}
saying the if $u\geq 0$, then
$$
H(u,Du,0,x)\leq \gb|Du|+\bar G,
$$
which yields \eqref{1.7.4}. Similarly \eqref{1.7.5}
is obtained. The lemma is proved. \qed

\begin{lemma}
                                          \label{lemma 1.6.2}
The functions $u_{n}$, $n=1,2,...$,
are uniformly continuous and uniformly bounded
in $\bar \Omega$ with the estimates of their sup
norms and moduluses of continuity involving
only the $L_{p}(\Omega)$-norm of $G$ and not its 
$L_{q}(\Omega)$-norm.
\end{lemma}

Proof. Fix $n$ and denote $\Omega_{+}=\{x\in\Omega:u_{n}(x)>0\}$.
By Lemma \ref{lemma 1.7.2} we have \eqref{1.7.4}
on $\Omega'$. By Theorem \ref{theorem 1.14.1} we have 
$$
u_{n}\leq N\|\bar G\|_{L_{p}(\Omega')}+
\sup_{\partial \Omega'}(u_{n})_{+}
\leq N\|\bar G\|_{L_{p}(\Omega)}+
\sup_{\partial \Omega'\cap \partial \Omega}g_{+},
$$
where $N$ depends only on $p,d,\delta, \|b\|$,
and the diameter of $D$. Similarly one estimates $u_{n}$
from below.

To show that $u_{n}$ are equicontinuous we use
\eqref{1.7.3}, denote $a^{ij}D_{ij}u_{n}+b^{i}D_{i}u_{n}=f_{n}$
and use that, in light of the first assertion,
 the $L_{p}(\Omega)$-norms of $f_{n}$
are uniformly bounded. Then by Theorem \ref{theorem 1.6.4}
  there is a constant $N$ depending only
on $p,d,\|b\|$, $\delta$, $\|\bar G\|_{L_{p}(\Omega)}$, $\sup|g|$,
and the diameter of $\Omega$, such that
$$
|u_{n}(x)-u_{n}(y)|\leq N(\rho(x)\wedge\rho(y))^{-\alpha}
|x-y|^{\alpha},
$$
where $\rho(z)$ is the distance from $z$ to $\Omega^{c}$
and $\alpha=\alpha(d,\delta,\|b\|)\in(0,1)$.

Furthermore, by Theorem \ref{theorem 1.6.1}
there is $\beta=\beta(d,\delta,\|b\|)>0$
and a constant $N$, depending only
on $ d,\delta,\|b\|  $, and  
  $\Omega$, such that
$$
|u_{n}(x)-u_{n}(x_{0})|\leq N|x-x_{0}|^{\beta}\|f_{n}\|_{L_{p}(\Omega)}
+w\big(N|x-x_{0}|^{\beta/2})
$$
whenever $x\in \Omega$ and $x_{0}\in \partial \Omega$,
where $w$ is the concave modulus of continuity of $g$.

A standard combination of these interior and boundary
estimates leads to our assertion. The lemma is proved.\qed

The main tool allowing us to pass to the limit under the sign of $H$
is given by the following lemma, which is stated
for the signs $\pm$ meaning that it holds when one takes
everywhere the upper sign and ignores the lower one
and also  holds when one takes everywhere the lower
sign and ignores the upper one.
It is worth saying that  generally, the results
of such kind are taken from Section 3.5 of \cite{Kr_85}. They    
generalize earlier results 
 for elliptic equations  by the author \cite{Kr_71}  
(1971) 
and
  Evans \cite{Ev_78}  (1978).
The methods in \cite{Kr_71} are quite transparent and are based on 
expressing the solution of the equation $H[u]=-f$
in the form $u=R_{\lambda}(\lambda u+f)$, $\lambda>0$,
where $R_{\lambda}$ is a nonlinear integral operator
continuous in $L_{p}$. It is easy to pass to the limit
under the sign of $R_{\lambda}$. In addition, it turns out that
if $u\in W^{2}_{p}$, then  $\lambda\big[R_{\lambda}(\lambda  u+f)-u]
\to F[u]+f$ as $\lambda\to
\infty$. Later on it became clear that the 
above integral representations
are equivalent to having \eqref{10.14.4} that possesses the same 
features
as the integral representations.

\begin{lemma}
                                          \label{lemma 1.7.10}
 Let $p\in[d_{0},\infty)$, $R\in (0,\infty)$,
$u\in W^{2}_{p }(B_{R})$.  Suppose that we are given
a function $H(\sfu'',x)$ that satisfies Assumption
\ref{assumption 3.11.2} and is such that
$$
|H(\sfu'',x)|\leq N_{0}|\sfu''|+\hat H(x),
$$
where $N_{0}$ is independent of $\sfu''$ and $x$
and nonnegative $\hat H$ belongs to $L_{p}(B_{R})$.
Then 

(i) there exists 
a constant
$N=N(p,d,\delta )\geq0$ such that 
for any $\lambda>0$ and $\phi\in  W^{2}_{p }(B_{R})$
we have
$$
\lambda\|(\phi-u)_{\pm}\|_{L_{p}(B_{R/2})}
\leq N\big\|\big(\lambda ( \phi-u)-(H [\phi]-f)\big)_{\pm}\big\|_{L_{p}(B_{R})}
$$
\begin{equation}
                                           \label{10.14.4}
+N\lambda R^{d/p}e^{-R\sqrt{\lambda}/N}\sup_{\partial B_{R}}
( \phi-u)_{\pm},
\end{equation}
for any $f\in L_{p}(B_{R})$ such that $(H [u]-f)_{\pm}=0$ on $B_{R}$;

(ii) if there is a constant $N$ such that \eqref{10.14.4}
holds for  an $f\in L_{p}(B_{R})$ and any sufficiently large
$\lambda>0$ and any $\phi\in  W^{2}_{p }(B_{R})$,
then $(H [u]-f)_{\pm}=0$  on $B_{R/2}$. 
\end{lemma}

Proof. (i) Observe that for an $\bS_{\delta}$-valued $a$
we have $H [\phi]-H[u]=a^{ij}D_{ij}(\phi-u)$.
Then the first assertion follows immediately from
Theorem \ref{theorem 10.14.2}.

(ii) Plug $u+  \phi/\lambda$ in \eqref{10.14.4}  in place of $\phi$.
Then
$$
 \| \phi _{+}\|_{L_{p}(B_{R/2})}
\leq N\big\|\big( \phi -( H [u+  
\phi/\lambda]-f)\big)_{+}\big\|_{L_{p}(B_{R})}
$$
$$
+N  R^{d/p}e^{-R\sqrt{\lambda}/N}\sup_{\partial B_{R}}
 \phi _{+}.
$$
Letting $\lambda\to\infty$ and using the dominated convergence theorem
 yields
$$
 \| \phi _{+}\|_{L_{p}(B_{R/2})}
\leq N\big\|\big( \phi -( H [u]-f)\big)_{+}\big\|_{L_{p}(B_{R})}.
$$
This is true for any $\phi\in  W^{2}_{p }(B_{R})$ and by continuity
for   any $\phi\in  L_{p }(B_{R})$.
Taking $\phi =  H [u]-f $ shows that $(H [u]-f)_{+}=0$.
This proves (ii) with the sign $+$. Similar argument is valid
for $-$. The lemma is proved.

{\bf Proof of Theorem \ref{theorem 10.5.1}}.
{\em Case $\bar G\in L_{q}(\Omega)$\/}.
Take the functions $u_{n}$ from Lemma \ref{lemma 1.6.1}
and extract a subsequence $u_{n(k)}$
such that 

(i) it converges weakly in $W^{2}_{p,\loc}(\Omega)$
 to a $u\in W^{2}_{p,\loc}(\Omega)$,
which is possible in light of Lemma \ref{lemma 1.6.1};

(ii) converges uniformly on $\Omega$ to $u$
thus making it belong to $C(\bar \Omega)$,
which is possible due to Lemma \ref{lemma 1.6.2};

(iii) there is a function $\hat H\in L_{d_{0},\loc}(\Omega)$ such 
that in $\Omega$
(a.e.)
\begin{equation}
                                                  \label{1.7.50}
\lim_{k\to\infty}H^{n(k)}_{u_{n(k)} ,Du_{n(k)} }( D^{2}u(x),x)
=H[u](x),
\end{equation}
\begin{equation}
                                                  \label{1.7.6}
\sup_{k}|H^{n(k)}_{u_{n(k)} ,Du_{n(k)} }( \sfu'',x)|
\leq N|\sfu''|+\hat H,
\end{equation}
which is possible due to Lemma \ref{lemma 1.7.1}.

Then for $m=1,2,...$ introduce
$$
\hat H^{m}(\sfu'',x)=\sup_{k\geq m}
H^{n(k)}_{u_{n(k)} ,Du_{n(k)} }( \sfu'',x).
$$
Obviously, for $k\geq m$ we have $\hat H^{m}[u_{n(k)}]\geq0$.

Due to Assumption \ref{assumption 3.11.2}, this assumption
is also satisfied for $\hat H^{m}$. Also, thanks to 
\eqref{1.7.6}, $\hat H^{m}$
satisfies the assumption of Lemma \ref{lemma 1.7.10}
with $p=d_{0}$.
By that lemma (with $f=0$ and sign $-$)
$$
\lambda\|(\phi-u_{n(k)})_{-}\|_{L_{d_{0}}(B_{R/2}(x))}
\leq N\big\|\big(\lambda ( \phi-u_{n(k)})-\hat
H^{m} [\phi]\big)_{-}\big\|_{L_{d_{0}}(B_{R}(x))}
$$
$$
+N\lambda R^{d/p}e^{-R\sqrt{\lambda}/N}\sup_{\partial B_{R}(x)}
( \phi-u_{n(k)})_{-},
$$
whenever $\bar B_{R}(x)\in \Omega$, $\phi\in W^{2}_{p}(B_{R}(x))$,
and $\lambda>0$. Passing to the limit as $k\to \infty$
and then using Lemma \ref{lemma 1.7.10} again
and using the arbitrariness of $B_{R}(x)$, we obtain
$\hat H^{m}[u]\geq 0$ (a.e.) in $\Omega$. Letting $m\to\infty$ and using 
\eqref{1.7.50} yields $H [u]\geq 0$ (a.e.) in $\Omega$.

One gets that $H [u]\leq 0$ (a.e.) in $\Omega$ similarly
by considering
$$
\check H^{m}(\sfu'',x)=\inf_{k\geq m}
H^{n(k)}_{u_{n(k)} ,Du_{n(k)} }( \sfu'',x)
$$
and using Lemma \ref{lemma 1.7.10} with sign $+$.
Finally estimate \eqref{1.8.2} follows from what is said
in the proof of Lemma \ref{lemma 1.6.2}.
This finishes the proof in our particular case
in which Assumption \ref{assumption 1.6.1} was not used.  

{\em General case\/}. 
In the above proof for $n\geq n_{0}$
(see Assumption \ref{assumption 1.6.1}) we replace $H^{n}$ 
introduced by \eqref{1.8.10} with
$$
H^{n}(\sfu,x)=H(\sfu,x)I_{\bar G(x)\leq n}+I_{\bar G(x)> n}F(
\sfu'_{0},\sfu'' ,x),
$$
 keep $F(\sfu'_{0},\sfu'',x)$ unchanged, and set
$G^{n}=H^{n}-F$. Then Assumptions  \ref{assumption 3.11.2},
\ref{assump1},   and \ref{assumption 10.5.1}
are obviously satisfied for the new couple $(H^{n},F)$
with the new $\bar G^{n}(x)=\bar G(x)I_{\bar G(x)\leq n}$
which is bounded. Assumption \ref{assumption 3,7,3}
is also satisfied with that $\bar G^{n}$.

After that literally repeating the above proofs 
 with the new $H^{n}$ proves the theorem
also in the general case.\qed
 
{\bf Acknowledgment}. The author is very grateful
to Hongjie Dong and A.I. Nazarov who read
the first draft of the paper and pointed
out several glitches in it.


\begin{thebibliography}{mm}

\bibitem{Ca_95}
X. Cabr\'e, {\em
On the Alexandroff-Bakelman-Pucci estimate and the reversed
H\"older inequality for solutions of elliptic 
and parabolic equations\/},
Comm. Pure Appl. Math., Vol. 48 (1995), 539--570.

\bibitem{Caf89} L.A. Caffarelli, {\em Interior a priori estimates 
for solutions of fully non-linear equations\/},  Ann. of Math.,  
Vol. 130 (1989), 189--213.

\bibitem{Caf90} L.A. Caffarelli, {\em 
Interior $W^{2,p}$
 estimates for solutions of the Monge-Amp\`ere equation\/}, 
Ann. of Math. (2), Vol. 131 (1990), No. 1, 135--150.

\bibitem{CC_95} L.A. Caffarelli and X. Cabr\'e, 
``Fully nonlinear elliptic
equations'',
 American Mathematical Society, Providence, 1995.

\bibitem{CCKS_96} L. Caffarelli, M. G. Crandall,
M. Kocan, and A. \'Swi{\c e}ch,,
{\em On viscosity
solutions of fully nonlinear equations with
 measurable ingredients\/},
 Comm.
Pure Appl. Math., Vol. 49 (1996), No. 4, 365--397.

\bibitem{Es93} L. Escauriaza, {\em $W^{2,n}$ a priori estimates for
solutions to fully non-linear equations\/}, Indiana Univ. Math. J.,
Vol. 42 (1993), No. 2, 413--423.

\bibitem{Ev_78} L.C. Evans,
{\em A convergence theorem for solutions of nonlinear
second order elliptic equations\/}, Indiana
University Math. J., Vol. 27 (1978), 875--887.

\bibitem{Fo_98} K. Fok, {\em A nonlinear Fabes-Stroock result\/},
Comm. PDEs, Vol 23 (1998), No. 5-6, 967--983.

\bibitem{KK_19}
Byungsoo Kang and Hyunseok Kim,  
{\em On $L^p$-resolvent estimates for second-order
elliptic equations in divergence form\/},
 Potential Anal., Vol. 50 (2019), No. 1,
107--133.

\bibitem{Kr_71}  N.V. Krylov, {\em
   On uniqueness of the solution of Bellman's equation\/}, 
Izvestiya Akademii Nauk SSSR, seriya matematicheskaya,
Vol. 35 (1971),  No. 6,  1377--1388 in Russian; English translation
in Math. USSR
Izvestija,
   Vol. 5 (1971), No. 6, 1387--1398.

\bibitem{Kr_85}  N.V. Krylov,
  ``Nonlinear elliptic and parabolic equations of second
  order'', Nauka, Moscow, 1985 in Russian;
English translation:  Reidel, Dordrecht, 1987. 

\bibitem{Kr_13.1} N.V. Krylov, {\em On the existence of $W^{2}_{p}$ 
solutions for fully nonlinear elliptic
equations under  relaxed convexity assumptions\/}, 
Comm. Partial Differential Equations, Vol. 38 (2013), No. 4,  
687--710.

\bibitem{Kr_18} N.V. Krylov,
``Sobolev and viscosity solutions for fully nonlinear  elliptic 
and parabolic equations'', Mathematical Surveys and Monographs,
233, Amer.
Math. Soc., Providence, RI, 2018.

 
\bibitem{Kr_19_1} N.V. Krylov,
{\em  On stochastic equations with drift in $L_{d}$\/}, 
\\ http://arxiv.org/abs/2001.04008


\bibitem{Kr_20} N.V. Krylov,
{\em On diffusion processes with drift in $L_{d}$\/},\\
http://arxiv.org/abs/2001.04950

\bibitem{KS_80}  N.V. Krylov and  M.V. Safonov,  {\em
   A certain property of solutions of parabolic equations
 with measurable coefficients}, 
Izvestiya Akademii Nauk SSSR, seriya matematicheskaya,
Vol.
44  (1980),  No. 1,  161--175  in Russian; English translation in
 Math. USSR
Izvestija, Vol. 16  (1981), No. 1,  151--164.

\bibitem{Sa_80}  M. V. Safonov, {\em
Harnack inequalities for elliptic equations and H{\"o}lder continuity
of their solutions\/}, 
Zap. Nauchn. Sem. Leningrad. Otdel. Mat. Inst. Steklov (LOMI),
Vol. 96 (1980), 272--287 in Russian;
English translation in Journal of Soviet Mathematics,
  Vol. 21 (March 1983), No. 5, 851--863.

  \bibitem{Sa_84}  M. V. Safonov, {\em
On the classical solutions of Bellman's elliptic
equations\/}, Dokl. Akad. Nauk SSSR, Vol. 278 (1984),
810--813 in Russian; English translation in Soviet Math. Dokl.,
Vol. 30 (1984), No. 2, 482--485.

\bibitem{Sa_88}  M. V. Safonov, {\em
On the classical solutions of nonlinear elliptic
equations of second order\/}, Izvestija Acad. Nauk SSSR, ser. matemat.,
Vol. 52 (1988), No. 6, 1272--1287   in Russian; English translation
in Math. USSR Izvestiya, Vol. 33 (1989), No. 3, 597--612.

\bibitem{Sa_10} M.V. Safonov, {\em Non-divergence elliptic
 equations of second
order with unbounded drift\/}, Nonlinear partial differential equations and
related topics, 211--232, Amer. Math. Soc. Transl. Ser. 2, 229, 
Adv. Math. Sci.,
64, Amer. Math. Soc., Providence, RI, 2010.
 
 \bibitem{SW_2016}
J. Streets and M. Warren, {\em Evans-Krylov estimates for a
nonconvex Monge-Amp\`ere equation\/}, Math. Ann., Vol. 
365 (2016), No. 1-2,
805--834.  

\bibitem{Wi09} N. Winter, {\em $W^{2,p}$ and $W^{1,p}$-estimates
at the boundary for solutions of fully nonlinear, uniformly elliptic equations\/},
 Z. Anal. Anwend., Vol. 28  (2009), No. 2, 129--164.


\end{thebibliography}
\end{document}